\documentclass{amsart}
\usepackage[T2A]{fontenc}
\usepackage[russian,english]{babel}
\def\a{\alpha}

\newtheorem{теор}{Тheorem}[section]
\newtheorem{лем}[теор]{Lemma}
\newtheorem{зам}[теор]{Remark}

\newtheorem{анн}[теор]{Annotation}
\newtheorem{опр}[теор]{Definition}
\newtheorem{ключ}[теор]{ Keywords.}

\newtheorem{thm}{Theorem}[section]

\begin{document}

\begin{center}
\textbf{{\large {\ 
Space-Time-Dependent Source Identification Problem for a Subdiffusion Equation}}}\\[0pt]
\medskip \textbf{R.R. Ashurov$^{1}$ and O.T. Mukhiddinova$^{1,2}$}\\[0pt]
\textit{ashurovr@gmail.com, oqila1992@mail.ru \\[0pt]}

\smallskip
\textit{{$^{1}$V.I. Romanovskiy Institute of Mathematics, Uzbekistan Academy of Science, University str.,9, Olmazor district, Tashkent, 100174, Uzbekistan}\\
{$^{2}$ Tashkent University of Information Technologies named after Muhammad al-Khwarizmi,
Str., 108, Amir Temur Avenue,Tashkent, 100200, Uzbekistan } }\\

\end{center}

\begin{анн}
In this paper, we investigate the inverse problem of determining the right-hand side of a subdiffusion equation with a Caputo time derivative, where the right-hand side depends on both time and certain spatial variables. Similar inverse problems have been previously explored for hyperbolic and parabolic equations, with some studies establishing the existence and uniqueness of generalized solutions, while others proved the uniqueness of classical solutions. However, such inverse problems for fractional-order equations have not been addressed prior to this work. Here, we establish the existence and uniqueness of the weak solution to the considered inverse problem. To solve it, we employ the Fourier method with respect to the variable independent of the unknown right-hand side, followed by the method of successive approximations to compute the Fourier coefficients of the solution. Notably, the results obtained are also novel for parabolic equations.
\end{анн}
\begin{ключ}
Subdiffusion equation, Caputo time derivative, inverse problem, uniqueness and existence of solution, Fourier method.
\end{ключ}
\section{Introduction}

The theory of differential equations with fractional derivatives has gained significant popularity and relevance in recent decades, largely due to their wide application across various, seemingly unrelated fields of science and technology (see, e.g., \cite{Machado}–\cite{SU}). Concurrently, the mathematical properties of fractional differential equations and methods for their solution have been thoroughly investigated by numerous researchers (see, e.g., \cite{Kil}–\cite{ACT}).

In recent years, there has been growing interest among researchers in inverse problems for differential equations of both integer and fractional order (see, e.g., the monograph by S. I. Kabanikhin \cite{Kaban} and the review article by Yamamoto \cite{Yamamoto1}). This interest in inverse problems stems from their critical role in applications across diverse fields such as mechanics, seismology, medical tomography, and geophysics (see, e.g., \cite{Machado}, \cite{Kaban}–\cite{Prilepko1}). 
Particular attention has been given to inverse problems of determining the source function when it takes the form \( F(x,t) = g(t)f(x) \), with either \( g(t) \) or \( f(x) \) unknown. To our knowledge, the more general case of \( F(x,t) \) without such a factorization has not been studied. In this scenario, even the choice of an appropriate overdetermination condition remains unclear. Inverse problems aimed at recovering the time-dependent component \( g(t) \) are typically addressed by reducing them to integral equations (see, e.g., \cite{ASh1}–\cite{ASh2} and references therein). In contrast, the inverse problem of determining the spatially dependent component \( f(x) \) is analyzed under two distinct cases: \( g(t) \equiv 1 \) and \( g(t) \not\equiv 1 \). When \( g(t) \equiv 1 \), such problems have been investigated, for example, in \cite{KirM}–\cite{AMux1}. The case where \( g(t) \not\equiv 1 \) is more intricate, with solvability depending on the properties of \( g(t) \). To tackle these problems, researchers commonly employ one of two overdetermination conditions: either \( u(x, T) = \psi(x) \) or an integral condition derived from the solution of the initial-boundary value problem (see, e.g., \cite{Prilepko1}, \cite{ASh3}–\cite{VanBockstal}, as well as \cite{Kaban} and the review article \cite{Yamamoto1}).

Inverse problems concerning the recovery of the right-hand side for equations of mixed type have been systematically studied by K. B. Sabitov and his students (see, e.g., \cite{Sabitov1}–\cite{Sabitov4}). These works explore the existence and uniqueness of the classical solution to such inverse problems.

The works most closely related to our study are those by S. Z. Dzhamalov et al. \cite{Djamalov1}–\cite{Djamalov5}, which address inverse problems of determining the right-hand side of equations, where the right-hand side depends on both time and a portion of the spatial variable. These studies investigate the existence of a generalized solution to the problem using the Galerkin method. Specifically, \cite{Djamalov1}–\cite{Djamalov3} consider various equations of mixed type, \cite{Djamalov4} examines a hyperbolic equation, and \cite{Djamalov5} focuses on a parabolic equation.

In \cite{Fikret}, the inverse problem for a hyperbolic equation is also addressed, specifically the determination of the right-hand side, which depends on both time and a portion of the spatial variable. The authors successfully established the uniqueness of the classical solution to this inverse problem.

Additionally, the work \cite{XiaomaoDeng} is noteworthy, where the authors proposed a numerical method for solving a similar problem for a parabolic equation.

In this paper, we investigate the inverse problem of determining the function \( h(t, x) \) that appears in the right-hand side of a subdiffusion equation, given as \( f(t, x, y) h(t, x) + g(t, x, y) \), where \( x \in \mathbb{R}^n \) and \( y \in \mathbb{R}^m \). For simplicity, we assume in the following that \( n = m = 1 \); the generalization to arbitrary \( n \) and \( m \) does not introduce additional fundamental difficulties.

Next, we proceed to the precise formulation of the problem.

 Let $ Q=(0,T) \times (0,1)\times (0,\pi)$
and
$ \Omega=(0,1)\times (0,\pi)$.
Consider the following initial - boundary value problem
\begin{equation} \label{1}
	\begin{cases}
		 D_t^{\alpha}u - u_{xx} - u_{yy} = g(t,x,y) + f(t,x,y) \cdot h(t,x) , \quad (t,x,y) \in Q, \\
        
         u(0,x,y) = \varphi(x,y), \quad (x,y) \in \Omega, \\
         
         u(t,0,y) =  u(t,1,y) = 0,\quad t\in [0,T], \quad y\in [0, \pi],\\
         
         u(t,x,0) = u(t,x,\pi) = 0,\quad t\in [0,T], \quad x\in [0,1].
	\end{cases}
\end{equation} 
 Here \( f, g \) and \( \varphi \) are given functions, \( \alpha \in (0,1) \), and \( D_t^{\alpha} \) is the fractional Caputo derivative, which for an absolutely continuous function \( v(t) \) is determined by the formula (see, e.g., \cite{Kil}, p. 91).
\[
D_{t}^{\alpha } v(t)\equiv J_t^{1-\alpha} D_t v(t), \quad J_t^\a v(t) = \frac{1}{\Gamma (\alpha )} \int _{0}^{t} (t-\tau )^{\alpha-1 }v(\tau ) d\tau ,
\] 
where $D_t=d/dt$, $\Gamma(\alpha)$ is the gamma function, $J_t^\alpha$ is the Riemann-Liouville fractional integral.

If the function \(h(t, x)\) is given, then the solution to the initial-boundary value problem \eqref{1} exists and is unique under certain conditions on the problem’s data (see, e.g., \cite{Sakamoto}, \cite{AMux2}).

Now suppose that the function \( h(t, x) \) is unknown and must be determined. The primary objective of this paper is to investigate the inverse problem of determining the pair of functions \( \{ u(t, x, y), h(t, x) \} \) subject to the overdetermination condition given by:
\begin{equation} \label{5}
u(t, x, l_0) = \psi(t, x), \quad t \in (0, T), \quad x \in (0, 1), \quad l_0 \in (0, \pi),
\end{equation}
where \( \psi \) is a known function.

This paper is organized into six sections. The following section contains supporting material, including well-known results of A. Alikhanov. Section 3 defines a weak solution to the inverse problem under consideration and presents the main result of the study. The solution to the problem is sought using the Fourier method. In Section 4, an infinite system of integro-differential equations is derived for the unknown Fourier coefficients. The method of successive approximations is employed to solve this system, and a priori estimates are established in this section. Section 5 proves the convergence of the corresponding sequences. Section 6 is dedicated to proving the main result. Finally, Section 7 provides the conclusion.

\section{Preliminaries}
In this section, we remind the definition of the Mittag - Leffler functions and introduce some auxiliary lemmas that will be used throughout the paper.

 The two-parameter Mittag-Leffler function $E_{\rho,\mu}(z)$ is an entire function defined by a power series of the form:
$$
E_{\rho,\mu}(z)= \sum\limits_{k=0}^\infty \frac{z^k}{\Gamma(\rho
k+\mu)},\quad \rho>0, \quad \mu, z\in \mathbb{C}.
$$
If \(\mu = 1\), then the Mittag-Leffler function is called the one-parameter or classical Mittag-Leffler function and is denoted by \(E_{\rho}(z) = E_{\rho,1}(z)\). Obviously, there is a constant $M$ such that
\begin{equation}\label{E}
   E_{\rho,\mu}(z) \leq M,  \quad \mu>0, \quad z\in [a,b]\subset \mathbb{R}_+. 
\end{equation}

We also need the following connection between the fractional integral and the derivative, which follows directly from the definitions and the equality for $v\in L_1(0,T)$: $J^\alpha_t(J^\beta_t v(t))= J^\beta_t(J^\alpha_t v(t))=J^{\alpha+\beta}_t v(t))$, $\alpha, \beta >0$ (see \cite{Dzh66}, p. 567).
\begin{лем}\label{J} Let $0<\alpha<1$ and $v(t)$ be absolutely continuous on $[0,T]$: $v\in AC\,[0,T]$. Then
\begin{equation}\label{JD}
    J_t^\alpha D^\alpha_t v(t)= v(t)- v(0).
\end{equation}
    
\end{лем}
Proof see in \cite{Dzh66}, p. 570.

Let us now recall the following two statements from the well-known work of Alikhanov \cite{Alixan}.
\begin{лем}\label{Alixan1}
Let $0<\alpha<1$ and  $w\in AC\, ([0,T]: L_2(0,1))$. Then
\[
\frac{1}{2} D_t^\alpha \| w(t,x) \|_{L_2(0, 1)}^2 \leq \int_0^1 w(t,x) \, D_t^\alpha w(t,x) \, dx.
\]
\end{лем}

\begin{лем}\label{Alixan2L}
Let \( y(t)\in AC\, [0,T] \) be a positive function and for all \( t \in (0,T] \), the following inequality holds:
\[
D_t^\alpha y(t) \leq c_1 y(t) + c_2(t), \quad 0<\alpha\leq 1,
\]
for almost all $t\in [0,T]$, where \( c_1>0 \) and \( c_2(t) \) is an integrable nonnegative function on $[0,T]$. Then
\begin{equation}\label{Alixan2}
y(t)\leq y(0) \,E_\alpha(c_1 t^\alpha) + \Gamma(\alpha) \,E_{\alpha, \alpha}(c_1 t^\alpha) \, J_t^\alpha c_2(t).
    \end{equation}
\end{лем}

Let us also note the following simple property of the fractional integral: Let $v(t)$ be a positive bounded function on $[0, T]$ and $\alpha\in (0,1]$. Then
\begin{equation}\label{integral}
\frac{1}{T^{1-\alpha}}\int_0^tv(s) ds \leq \int_0^t\frac{v(s)}{t^{1-\alpha}} ds\leq \int_0^t\frac{v(s)}{(t-s)^{1-\alpha}}ds=\Gamma(\alpha) J_t^\alpha v(t)\leq T^\alpha \max\limits_{t\in [0,T]} v(t).
    \end{equation}

\section{Definition of the Weak Solution and Formulation of the Main Result}
The solution to the inverse problem \eqref{1}–\eqref{5} is sought in the form of a formal Fourier series:
\begin{equation}\label{6}
    u(t, x, y) = \sum_{k=1}^\infty u_k(t, x) \sin k y,
\end{equation}
where \( u_k(t, x) \) are the unknown functions and \( \sin k y \) are the eigenfunctions of the spectral problem:
\[
    -Y''(y) + \lambda Y(y) = 0, \quad
    Y(0) = Y(\pi) = 0.
\]
The corresponding eigenvalues are given by \( \lambda_k = k^2 \). Here after, the symbol \( v_k \) denotes the Fourier coefficients of a function \( v(y) \in L_1(0, \pi) \) with respect to the system \( \{\sin k y\} \), defined as:
\begin{equation}\label{CoefFur}
    v_k = \frac{2}{\pi} \int_0^\pi v(y) \sin k y \, dy.
\end{equation}

Next, we apply the method proposed by B. A. Bubnov \cite{Bub} and further developed by S. Z. Dzhamalov and other authors (see \cite{Djamalov1}–\cite{Djamalov5}). According to this method, we substitute \( y = l_0 \) into equation \eqref{1}, yielding:
\begin{equation*}
    D_t^{\alpha} u(t, x, l_0) - u_{xx}(t, x, l_0) - u_{yy}(t, x, l_0) = f(t, x, l_0) h(t, x) + g(t, x, l_0).
\end{equation*}
Using the form of the function \( u(t, x, y) \) given by \eqref{6} and applying the overdetermination condition \eqref{5}, we obtain:
\begin{equation*}
    D_t^{\alpha} \psi(t, x) - \psi_{xx}(t, x) - \sum_{k=1}^{\infty} \lambda_k u_k(t, x) \sin k l_0 = f(t, x, l_0) h(t, x) + g(t, x, l_0).
\end{equation*}
Assuming that \( f(t, x, l_0) \neq 0 \) for all \( (t, x) \in [0, T] \times [0, 1] \), we derive a formal expression for \( h(t, x) \):
\begin{equation}\label{7}
    h(t, x) = \frac{D_t^{\alpha} \psi(t, x) - \psi_{xx}(t, x)-g(t, x, l_0) - \sum_{k=1}^{\infty} \lambda_k u_k(t, x) \sin k l_0}{f(t, x, l_0)}.
\end{equation}

Let \( B \) be a Banach space. We denote by \( L_\infty(0, T; B) \) the space of functions that are essentially bounded on \( (0, T) \) and take values in \( B \). The space \( L_1(0, T; B) \) is defined similarly. Let \( W_2^k(\Omega) \) denote a classical Sobolev space. Then, the symbol \( \dot{W}_2^1(\Omega) \) represents the closure of the set \( C_0^\infty(\Omega) \) with respect to the norm of \( W_2^1(\Omega) \). Hereafter, \( \|\cdot\|_{_2} \) denotes the norm in \( L_2(\Omega) \).

We now define a weak formulation of the inverse problem \eqref{1}–\eqref{5}.

\begin{опр}\label{opr1} 
Find a pair of functions \( \{ u(t, x, y), h(t, x) \} \), where \( h(t, x) \), with Fourier coefficients \( u_k(t, x) \) determined by formula \eqref{CoefFur}, takes the form given in \eqref{7}, and the function \( u(t, x, y) \) satisfies the following conditions:
\begin{enumerate}
    \item \( u \in L_\infty(0, T; L_2(\Omega)) \), \,\, \( u \in L_1(0, T; \dot{W}_2^1(\Omega)) \);
    \item \( D_t^\alpha u \in L_1(0, T; L_2(\Omega)) \);
    \item \( u(0, x, y) = \varphi(x, y) \) a.e.\ in \( \Omega \);
    \item For any \( v \in \dot{W}_2^1(\Omega) \) and almost every \( t \in (0, T] \), the following equality holds:
\end{enumerate}
\begin{equation}\label{equation}
    \int_\Omega D_t^\alpha u v \, dx \, dy + \int_\Omega (u_x v_x + u_y v_y) \, dx \, dy = \int_\Omega f h v \, dx \, dy + \int_\Omega g v \, dx \, dy.
\end{equation}
\end{опр}

Let $M_\alpha = M T^\alpha$, where $M$ is from (\ref{E}). We now state the main result of the paper.

\begin{thm}\label{theorem 1}

Suppose the following conditions hold:
\begin{enumerate}
    \item \( f(t, x, l_0) \in C([0, T] \times [0, 1]) \), \( f(t, x, l_0) \neq 0 \), and \( f_0 = \max_{t, x} \left| \frac{1}{f(t, x, l_0)} \right| \);
    \item \( g(t, x, l_0) \in C([0, T] \times [0, 1]) \), and \( g_0 = \max_{t, x} |g(t, x, l_0)| \);
    \item \( D_t^\alpha \psi, \psi_{xx} \in C([0, T] \times [0, 1]) \), and \( \psi_0 = \max_{t, x} \left( |D_t^\alpha \psi| + |\psi_{xx}| \right) \);
    \item   \( 2 M_\alpha f_0^2 \max_{t,x}\| f_{yyy}(t, x, \cdot) \|_{L_2(0, \pi)}^2 \leq 1 \) and \(f(t, x, 0) = f(t, x, \pi) =  f_{yy}(t, x, 0) = f_{yy}(t, x, \pi) = 0 \);
    \item \( \max_{t}\| g_{yyy}(t, \cdot, \cdot) \|_{2} < \infty \), \, and \(g(t, x, 0) = g(t, x, \pi) =  g_{yy}(t, x, 0) = g_{yy}(t, x, \pi) = 0 \);
     \item \( \| \varphi_{x} \|_{2},\,\| \varphi_{yyy} \|_{2} < \infty, \) and \,\(\varphi(t, x, 0) = \varphi(t, x, \pi) =  \varphi_{yy}(t, x, 0) = \varphi_{yy}(t, x, \pi) = 0 \).
\end{enumerate}
Then, the inverse problem has a unique weak solution. Moreover, the following estimates hold:
\begin{equation}\label{uestimate1}
\max_{t}\| u(t, \cdot, \cdot) \|_2^2\leq 2 B_1,
\end{equation}
 \begin{equation}\label{u_xu_y}
|| u_x ||_{L_2(Q)}^2 +|| u_y ||_{L_2(Q)}^2 \leq \frac{\Gamma(\alpha)T^{1-\alpha}}{2}||\varphi||^2_2+ 3 B_1 T +\frac{T}{2} \, A_1
\end{equation}
\[
+T\, B_1 f_0^2 \max_{t, x} ||f(t,x,\cdot)||_{L_2(0, 
\pi)}^2, 
\]
\begin{equation}\label{DalphaInt}
    || D_t^\alpha u||_{L_2(Q)}^2 \leq \frac{\Gamma(\alpha)T^{1-\alpha}}{2} \left(\| \varphi_x \|^2_2 +  \| \varphi_y \|^2_2\right)
\end{equation}
\[
+T\max_t\left[f_0 ^2(\psi_0 + g_0)^2 \|f(t, \cdot, \cdot)\|^2_2
+ \|g(t, \cdot, \cdot)\|^2_2\right]+2 B_1 f_0^2\max_{t, x} ||f(t,x,\cdot)||_{L_2(0, \pi)}^2 , 
\]

\begin{equation}\label{h}
   \max_{t} ||h(t,\cdot)||_{L_2(0,1)}^2 \leq 4 f_0^2\left(\psi_0^2 +  2 B_1\right)+ 2 g_0^2.
\end{equation}
Here
\[
A_1= f_0^2 (\psi_0 + g_0)^2 \max_{t}\|f(t, \cdot, \cdot)\|_{L_2(\Omega)}^2 + \max_{t}\|g(t, \cdot, \cdot)\|_{L_2(\Omega)}^2,
\]
\[
B_1 = M \| \varphi_{yyy} \|_{L_2(\Omega)}^2 + M_\alpha f_0^2(\psi_0 + g_0)^2  \max_{t}\| f_{yyy}(t, \cdot, \cdot) \|_{L_2(\Omega)}^2+ M_\alpha  \max_{t}\| g_{yyy}(t, \cdot, \cdot) \|_{L_2(\Omega)}^2.
\]
\end{thm}

\begin{зам}\label{zam.1} By virtue of S. M. Bernstein's theorem on the uniform convergence of trigonometric series (see \cite{Zyg}, p. 384), the conditions of the theorem on the functions $f, g, \varphi$ guarantee absolute and uniform convergence on \( [0, \pi] \) of the corresponding Fourier series.
\end{зам}

\begin{зам}\label{zam.2}In fact, the theorem is proven for a wider class of functions $f, g, \varphi$. Specifically, its validity is established under the following conditions for some \( \varepsilon > 0 \) and all \( t \in [0, T] \):
\[
M_\alpha\frac{f_0^2}{\varepsilon} \max_{t,x} \sum_{k=1}^{\infty} \lambda_k^{5/2 + \varepsilon} |f_k(t, x)|^2 \leq 1.
 \quad \text{(in place of condition (4) of Theorem \ref{theorem 1})},
\]
\[
\max_t\sum_{k=1}^\infty \lambda_k^{\frac{5}{2} + \varepsilon} \|g_k(t, \cdot)\|_{L_2(0, 1)}^2 < \infty \quad \text{(in place of condition (5) of Theorem \ref{theorem 1})},
\]
\[
\sum_{k=1}^\infty \lambda_k^{\frac{5}{2} + \varepsilon} \|\varphi_k\|_{L_2(0, 1)}^2 < \infty \quad \text{(in place of condition (6) of Theorem \ref{theorem 1})}.
\]
If these conditions are met, a stronger property of the solution than (\ref{uestimate1}) is proven, i.e. the estimate is  valid:
\begin{equation}\label{182}
\max_{[0,T]}\sum_{k=1}^\infty \lambda_k^{\frac{5}{2} + \varepsilon} \| u_k \|_{L_2(0, 1)}^2 \leq 2 A_0,
\end{equation}

If we put \( \varepsilon = \frac{1}{2} \) in the above conditions, then we obtain the conditions of the theorem.
\end{зам}

Indeed, suppose, for example, that condition (5) of the theorem holds and \( \varepsilon = \frac{1}{2} \). Then integrating by parts, we obtain:
\begin{align*}
g_k(t, x) &= \frac{2}{\pi} \int_0^\pi g(t, x, y) \sin k y \, dy = -\frac{2}{k \pi} \int_0^\pi \frac{\partial g}{\partial y} \cos k y \, dy \\
&= \frac{2}{k^2 \pi} \int_0^\pi \frac{\partial^2 g}{\partial y^2} \sin k y \, dy=-\frac{2}{k^3 \pi} \int_0^\pi \frac{\partial^3 g}{\partial y^3} \cos k y \, dy =: \frac{1}{k^2} \tilde{g}_{yyy,k}(t, x).
\end{align*}
Thus, recalling that \( \lambda_k = k^2 \) and \( \varepsilon = \frac{1}{2} \), it follows from  Bessel’s inequality that:
\[
\sum_{k=1}^\infty \lambda_k^3 \|g_k(t, \cdot)\|_{L_2(0, 1)}^2 = \int_0^1 \sum_{k=1}^\infty |\tilde{g}_{yyy,k}(t, x)|^2 \, dx \leq \frac{2}{\pi} \|g_{yyy}(t, \cdot, \cdot)\|_{2}^2.
\]

\section{A priori estimates}

We decompose the functions \( f(t, x, y) \), \( g(t, x, y) \), and \( \varphi(x, y) \) into Fourier series (see Remark \ref{zam.1}) and denote their corresponding Fourier coefficients by \( f_k(t, x) \), \( g_k(t, x) \), and \( \varphi_k(x) \), respectively (cf. \eqref{CoefFur}). Substituting these series and the representation \eqref{6} for \( u(t, x, y) \) into equation \eqref{equation}, with the test function \( v(x, y) \) replaced by \( w(x) \sin k y \), where \( w \in \dot{W}_2^1(0, 1) \), we obtain the following equation to determine the unknown coefficients \( u_k(t, x) \) (see \eqref{6}):
\begin{equation}\label{81}
\int_0^1 D_t^\alpha u_k w \, dx + \int_0^1 (u_k)_x w_x \, dx + \lambda_k \int_0^1 u_k w \, dx = \int_0^1 (g_k + f_k h) w \, dx,
\end{equation}
where \( h \) is given by \eqref{7}. The challenge in solving \eqref{81} for \( u_k \) arises because \( h \) depends on all \( u_j(t, x) \), \( j = 1, 2, \dots \). To address this, we employ the method of successive approximations by constructing the sequence \( \{ u_k^n \} \), \( n = 1, 2, \dots \), and then proving that \( \{ u_k^n \} \to u_k \) as \( n \to \infty \) in an appropriate norm.

Thus, using the explicit expression of the function \( h \) (see \eqref{7}), we derive an initial-boundary value problem in the form of a recurrence relation for all \( k \geq 1 \) and \( n \geq 1 \):
\begin{equation}\label{9}
\int_0^1 D_t^\alpha u_k^n w \, dx + \int_0^1 (u_k^n)_x w_x \, dx + \lambda_k \int_0^1 u_k^n w \, dx 
\end{equation}
\[
= \int_0^1 M_k(t, x) w \, dx - \int_0^1 \frac{f_k(t, x)}{f(t, x, l_0)} \sum_{j=1}^\infty \lambda_j u_j^{n-1}(t, x) \sin j l_0 w \, dx,
\]
for all \( w \in \dot{W}_2^1(0, 1) \), with the initial condition:
\begin{equation}\label{10}
u_k^n(0, x) = \varphi_k(x),
\end{equation}
where
\[
M_k(t, x) = f_k(t, x) \left[ \frac{D_t^\alpha \psi(t, x) - \psi_{xx}(t, x) - g(t, x, l_0)}{f(t, x, l_0)} \right] + g_k(t, x).
\]

Note that if the right-hand side of equation \eqref{9} is known, the problem \eqref{9}–\eqref{10} is well studied (see, e.g., \cite{Sakamoto}, \cite{AMux2}). In particular, \cite{Sakamoto} establishes that if the right-hand side (for all \( t \in (0, T] \)) and the initial function \( \varphi_k(x) \) belong to \( L_2(0, 1) \) (as is the case here), then there exists a unique strong solution to the problem such that \( u_k^n(t, x) \in L_\infty(0, T; \dot{W}_2^2(0, 1)) \) and \( D_t^\alpha u_k^n(t, x) \in L_\infty(0, T; L_2(0, 1)) \). Naturally, such a solution also satisfies the problem \eqref{9}–\eqref{10}.

We analyze the problem \eqref{9}–\eqref{10} as follows. We take \( u_k^0 = 0 \) as the initial approximation for all \( k = 1, 2, \dots \). Then, we solve the initial-boundary value problem \eqref{9}–\eqref{10} to obtain all \( \{ u_k^1 \} \), \( k = 1, 2, \dots \). Next, we construct the sequence \( \{ u_k^n \} \) iteratively. However, before doing so, assuming that all \( \{ u_k^n \} \) have been constructed, we first establish a priori estimates.

\begin{лем}\label{u_k} Let $\varepsilon>0$. If
\begin{equation}\label{fk}
M_\alpha\frac{f_0^2}{\varepsilon} \max_{t,x} \sum_{k=1}^{\infty} \lambda_k^{5/2 + \varepsilon} |f_k(t, x)|^2 \leq 1.
\end{equation}
then
\begin{equation}\label{181}
\max_{[0,T]}\sum_{k=1}^\infty \lambda_k^{\frac{5}{2} + \varepsilon} \| u_k^n \|_{L_2(0, 1)}^2 \leq 2 A_0,
\end{equation}
where 
\[
A_0 = M \varphi^* + M_\alpha f_0^2(\psi_0 + g_0)^2 f^*+ M_\alpha g^*, \,\,
\varphi^*= \int_0^1 \sum_{k=1}^{\infty} \lambda_k^{5/2 + \varepsilon} |\varphi_k|^2 dx, 
\]
and
\[
f^*=\max_{[0,T]}\int_0^1 \sum_{k=1}^{\infty} \lambda_k^{\frac{5}{2}+\varepsilon} |f_k(t,x)|^2 dx, \quad g^*=\max_{[0,T]}\int_0^1 \sum_{k=1}^{\infty} \lambda_k^{5/2 + \varepsilon} |g_k(t,x)|^2 dx.
\]

\end{лем}

\begin{proof} Substitute \( w = u_k^n \) into equation \eqref{9}, yielding:
\begin{equation}\label{82}
\int_0^1 D_t^\alpha u_k^n u_k^n \, dx + \| (u_k^n)_x \|_{L_2(0, 1)}^2 + \lambda_k \| u_k^n \|_{L_2(0, 1)}^2
\end{equation}
\[
 = \int_0^1 M_k(t, x) u_k^n \, dx - \int_0^1 u_k^n \frac{f_k(t, x)}{f(t, x, l_0)} \sum_{j=1}^\infty \lambda_j u_j^{n-1} \sin j l_0 \, dx.
\]
For the integral on the left-hand side, we apply Alikhanov’s estimate (see Lemma \ref{Alixan1}), obtaining:
\begin{equation}\label{12}
\frac{1}{2} D_t^\alpha \| u_k^n \|_{L_2(0, 1)}^2 + \| (u_k^n)_x \|_{L_2(0, 1)}^2 + \lambda_k \| u_k^n \|_{L_2(0, 1)}^2 
\end{equation}
\[
\leq \left| \int_0^1 M_k(t, x) u_k^n \, dx \right| + \left| \int_0^1 u_k^n \frac{f_k(t, x)}{f(t, x, l_0)} \sum_{j=1}^\infty \lambda_j u_j^{n-1} \sin j l_0 \, dx \right| := I_1(t) + I_2(t).
\]
We now estimate the right-hand side. Using the conditions and notation of Theorem \ref{theorem 1}, we have:
\[
I_1(t) \leq f_0 (\psi_0 + g_0) \int_0^1 |f_k(t, x) u_k^n(t, x)| \, dx + \int_0^1 |g_k(t, x) u_k^n(t, x)| \, dx.
\]
Applying the inequality \( 2ab \leq a^2 + b^2 \), this becomes:
\begin{equation}\label{13}
I_1(t) \leq \frac{1}{2}\left[f_0^2 (\psi_0 + g_0)^2 \| f_k \|_{L_2(0, 1)}^2 + \| g_k \|_{L_2(0, 1)}^2\right]+  \| u_k^n \|_{L_2(0, 1)}^2.
\end{equation}

Further, using the notation \( f_0 \) from Theorem \ref{theorem 1}, we estimate:
\[
I_2(t)\leq \frac{1}{2} \| u_k^n \|_{L_2(0, 1)}^2 +  \frac{1}{2}f_0^2 \int_0^1 |f_k(t, x)|^2 \left| \sum_{j=1}^\infty \lambda_j u_j^{n-1} \sin j l_0 \right|^2 \, dx.
\]
For any \( \varepsilon > 0 \), write \( \lambda_k = \lambda_k^{1 + \frac{1 + 2\varepsilon}{4}} \lambda_k^{-\frac{1 + 2\varepsilon}{4}} \) and apply the Cauchy–Bunyakovsky inequality to the sum, yielding:
\[
I_2(t)\leq \frac{1}{2} \| u_k^n \|_{L_2(0, 1)}^2 +  \frac{1}{2}f_0^2 \int_0^1 |f_k(t, x)|^2 \sum_{j=1}^\infty \lambda_j^{-\frac{1}{2} - \varepsilon} \sum_{j=1}^\infty \lambda_j^{\frac{5}{2} + \varepsilon} |u_j^{n-1}|^2 \, dx.
\]
Since \( \sum_{j=1}^\infty \lambda_j^{-\frac{1}{2} - \varepsilon} \leq \frac{1}{2\varepsilon} \), this becomes:
\begin{equation}\label{14}
I_2(t)\leq \frac{1}{2} \| u_k^n \|_{L_2(0, 1)}^2 + \frac{f_0^2}{4\varepsilon} \int_0^1 |f_k(t, x)|^2 \sum_{j=1}^\infty \lambda_j^{\frac{5}{2} + \varepsilon} |u_j^{n-1}|^2 \, dx.
\end{equation}

Combining \eqref{13} and \eqref{14}, the estimate in \eqref{12} can be rewritten as:
\begin{equation}\label{15}
\frac{1}{2} D_t^\alpha \| u_k^n \|_{L_2(0, 1)}^2 + \| (u_k^n)_x \|_{L_2(0, 1)}^2 + \lambda_k \| u_k^n \|_{L_2(0, 1)}^2\leq \frac{3}{2} \| u_k^n \|_{L_2(0, 1)}^2 
\end{equation}
\[
 +  \frac{1}{2}\left[f_0^2 (\psi_0 + g_0)^2 \|f_k(t, \cdot)\|_{L_2(0, 1)}^2 + \|g_k(t, \cdot)\|_{L_2(0, 1)}^2\right] + \frac{f_0^2}{4\varepsilon} \int_0^1 |f_k(t, x)|^2 \sum_{j=1}^\infty \lambda_j^{\frac{5}{2} + \varepsilon} |u_j^{n-1}|^2 \, dx.
\]
If we discard the last two terms on the left - hand side, then
\[
D_t^\alpha \| u_k^n \|_{L_2(0, 1)}^2 \leq 3\| u_k^n \|_{L_2(0, 1)}^2 +c^k_2(t)
\]
where $ c_2^k(t)= c^k_{2,1}(t)+ c^k_{2, 2}(t)$ and
\[
 c^k_{2,1}(t)= f_0^2 (\psi_0 + g_0)^2 \|f_k(t, \cdot)\|_{L_2(0, 1)}^2 + \|g_k(t, \cdot)\|_{L_2(0, 1)}^2,
 \]
 \[
c^k_{2, 2}(t)= \frac{f_0^2}{2\varepsilon} \int_0^1 |f_k(t, x)|^2 \sum_{j=1}^\infty \lambda_j^{\frac{5}{2} + \varepsilon} |u_j^{n-1}|^2 \, dx.
\]
Apply Lemma \ref{Alixan2L} to obtain
\[
\| u_k^n \|_{L_2(0, 1)}^2\leq \| \varphi_k \|_{L_2(0, 1)}^2 \,E_\alpha(3 t^\alpha) + \Gamma(\alpha) \,E_{\alpha, \alpha}(3 t^\alpha) \, J_t^\alpha c^k_2(t).
\]
Based on the boundedness of the Mittag-Leffler function (see (\ref{E})) and estimate (\ref{integral}), it follows that
\[
\| u_k^n(t, \cdot) \|_{L_2(0, 1)}^2\leq M\| \varphi_k \|_{L_2(0, 1)}^2  + MT^\alpha \left(\max_{[0,T]}c^k_{2,1}(t) +\max_{[0,T]}c^k_{2,2}(t)\right).
\]
If we multiply this inequality by $\lambda_k^{\frac{5}{2}+\varepsilon}$, and then sum over $k$ from 1 to $\infty$, we get
\begin{align*}
&\max_{[0,T]}\sum_{k=1}^{\infty} \lambda_k^{5/2 + \varepsilon} \| u_k^n (t, \cdot)\|^2_{L_2(0,1)} \leq  M\int_0^1 \sum_{k=1}^{\infty} \lambda_k^{5/2 + \varepsilon} |\varphi_k|^2 dx + M_\alpha f_0^2(\psi_0 + g_0)^2 f^*\\
&+ M_\alpha g^*+ M_\alpha\frac{f_0^2}{2\varepsilon} \max_{[0,T]}\int_0^1 \sum_{k=1}^{\infty} \lambda_k^{5/2 + \varepsilon} |f_k(t, x)|^2 \sum_{j=1}^{\infty} \lambda_j^{\frac{5}{2}+\varepsilon} | u_j^{n-1}(t, x)|^2 dx.
\end{align*}
Using condition (\ref{fk}) and the notation of Lemma \ref{u_k}
the last estimate can be rewritten in the form of a recurrence estimate:
\begin{equation}\label{17}
\max_{[0,T]}\sum_{k=1}^{\infty} \lambda_k^{\frac{5}{2} + \varepsilon} \| u_k^n \|_{L_2(0,1)}^2 \leq A_0 + \frac{1}{2} \max_{[0,T]}\sum_{k=1}^{\infty} \lambda_k^{\frac{5}{2} + \varepsilon} \| u_k^{n-1} \|_{L_2(0,1)}^2, \quad n = 1, 2, \dots.
\end{equation}

As noted above, we take \( u_k^0 = 0 \) as the initial approximation for all \( k \geq 1 \). Then, for \( u_k^1(t, x) \), \( k = 1, 2, \dots \), from \eqref{17}, we obtain:
\begin{equation}\label{u1}
\max_{[0,T]}\sum_{k=1}^\infty \lambda_k^{\frac{5}{2} + \varepsilon} \| u_k^1 \|_{L_2(0, 1)}^2 \leq A_0.
\end{equation}
Next, we substitute these functions \( \{ u_k^1 \}_{k=1}^\infty \) into the problem \eqref{9}–\eqref{10} to uniquely determine the functions \( \{ u_k^2 \}_{k=1}^\infty \). For these functions, we derive the following estimate from \eqref{17}:
\[
\max_{[0,T]}\sum_{k=1}^\infty \lambda_k^{\frac{5}{2} + \varepsilon} \| u_k^2 \|_{L_2(0, 1)}^2 \leq A_0 + \frac{1}{2} A_0.
\]
Continuing this process, we finally obtain:
\begin{equation}\label{18}
\max_{[0,T]}\sum_{k=1}^\infty \lambda_k^{\frac{5}{2} + \varepsilon} \| u_k^n \|_{L_2(0, 1)}^2 = A_0 \sum_{j=1}^n \left( \frac{1}{2} \right)^j = 2 \left( 1 - \frac{1}{2^n} \right) A_0 \leq 2 A_0.
\end{equation}
\end{proof}

\begin{лем}One has
\[
\sum\limits_{k=1}^\infty\int\limits_0^t\int\limits_0^1\left[| (u_k^n)_x |^2 + \lambda_k | u_k^n |^2\right]dx d\tau \leq \frac{\Gamma(\alpha)T^{1-\alpha}}{2}||\varphi||^2_2+ 3 A_0 T +\frac{T}{2} \, A_1
\]
\begin{equation}\label{u_x}
+\frac{T\, A_0 f_0^2}{2 \varepsilon} \max_{x\in [0,1]} \int\limits_0^\pi |f(t,x,y)|^2 dy, \quad t\in [0,T]. 
\end{equation}

\end{лем}

\begin{proof}

Apply the operator $J_t^\alpha$ to both parts of inequality (\ref{15}) and use Lemma \ref{J}. Then
\begin{equation}\label{151}
\| u_k^n \|_{L_2(0, 1)}^2 + 2 J_t^\alpha\| (u_k^n)_x \|_{L_2(0, 1)}^2 + 2 \lambda_k J_t^\alpha\| u_k^n \|_{L_2(0, 1)}^2 
\end{equation}
\[
\leq \| \varphi_k \|_{L_2(0, 1)}^2+3 J_t^\alpha \| u_k^n \|_{L_2(0, 1)}^2 +  J_t^\alpha c_{2,1}^k(t) +J_t^\alpha c_{2,2}^k(t).
\]
First, we omit the first term on the left-hand side of this inequality. Then, summing inequality \eqref{151} over \( k \) from 1 to \( \infty \) and applying the estimates (\ref{integral}), we obtain
\begin{equation}\label{151}
\frac{2}{\Gamma(\alpha)T^{1-\alpha}} \sum\limits_{k=1}^\infty\left[\int\limits_0^t\int\limits_0^1| (u_k^n)_x |^2dxd\tau +  \lambda_k \int\limits_0^t\int\limits_0^1| u_k^n |^2dxd\tau\right] 
\end{equation}
\[
\leq \sum\limits_{k=1}^\infty\| \varphi_k \|_{L_2(0, 1)}^2+\frac{3 T^\alpha}{\Gamma(\alpha)}\max_{[0,T]} \sum\limits_{k=1}^\infty\| u_k^n \|_{L_2(0, 1)}^2 +  \frac{ T^\alpha}{\Gamma(\alpha)}\max_{[0,T]} \sum\limits_{k=1}^\infty \left(c_{2,1}^k(t) + c_{2,2}^k(t)\right).
\]
From Parseval's equality, it follows
\[
\max_{[0,T]}\sum\limits_{k=1}^\infty c_{2,1}^k(t)=f_0^2 (\psi_0 + g_0)^2 \max_{[0,T]}\|f(t, \cdot, \cdot)\|_{L_2(\Omega)}^2 + \max_{[0,T]}\|g(t, \cdot, \cdot)\|_{L_2(\Omega)}^2:= A_1.
\]
Similarly,
\[
\sum\limits_{k=1}^\infty c_{2,2}^k(t) \leq \frac{f_0^2}{2\varepsilon} \max_{x\in [0,1]} \int\limits_0^\pi |f(t,x,y)|^2 dy \sum_{k=1}^{\infty} \lambda_k^{\frac{5}{2} + \varepsilon} \| u_k^{n-1}(t, \cdot) \|_{L_2(0,1)}^2,
\]
or by (\ref{181}),
\[
\sum\limits_{k=1}^\infty c_{2,2}^k(t) \leq  \frac{A_0 f_0^2}{\varepsilon} \max_{x\in [0,1]} \int\limits_0^\pi |f(t,x,y)|^2 dy. 
\]
Therefore, inequality (\ref{151}) can be rewritten as (\ref{u_x}).
\end{proof}

\begin{лем} One has
\begin{equation}\label{Dalpha1}
\sum\limits_{k=1}^\infty\int\limits_0^t\int\limits_0^1| D_t^\alpha u_k^n |^2dxd\tau \leq \frac{\Gamma(\alpha)T^{1-\alpha}}{2} \left(\| \varphi_x \|^2_2 +  \| \varphi_y \|^2_2\right)
\end{equation}
\[
+T\max_t\left[f_0 ^2(\psi_0 + g_0)^2 \|f(t, \cdot, \cdot)\|^2_2
+ \|g(t, \cdot, \cdot)\|^2_2\right]+\frac{A_0 f_0^2}{\varepsilon} \max_{x\in [0,1]} \int\limits_0^\pi |f(t,x,y)|^2 dy. 
\]
\end{лем}
\begin{proof}
Substitute \( w = D_t^\alpha u_k^n \) into equation \eqref{9}, yielding:
\begin{equation}\label{82}
\| D_t^\alpha u_k^n \|_{L_2(0, 1)}^2 + \int_0^1 (u_k^n)_x D_t^\alpha (u_k^n)_x \, dx + \lambda_k \int_0^1 u_k^n D_t^\alpha u_k^n \, dx 
\end{equation}
\[
= \int_0^1 M_k(t, x) D_t^\alpha u_k^n \, dx - \int_0^1 D_t^\alpha u_k^n \frac{f_k(t, x)}{f(t, x, l_0)} \sum_{j=1}^\infty \lambda_j u_j^{n-1} \sin j l_0 \, dx.
\]
Let us apply A. Alikhanov's estimate (see, Lemma \ref{Alixan1}) in the second and third integrals on the left-hand side of the equality. Then
\begin{equation}\label{121}
   \| D_t^{\alpha} u_k^n \|^2_{L_2(0,1)} + \frac{1}{2} D_t^{\alpha} \| (u_k^n)_x \|^2_{L_2(0,1)} + \lambda_k \frac{1}{2} D_t^{\alpha} \| u_k^n \|^2_{L_2(0,1)}  
\end{equation}
\[
\leq \left|\int_0^1 M_k(t,x) D_t^{\alpha} u_k^n dx\right| + \left|\int_0^1 D_t^{\alpha} u_k^n \frac{f_k (t,x)}{f(t,x,l_0)} \sum_{j=1}^{\infty} \lambda_j u_j^{n-1} \sin jl_0 dx\right| := J_1(t) + J_2(t).
\]
For $J_1(t)$ one has (see Eq. \eqref{13})
\[
J_1(t)  \leq  f_0^2(\psi_0  + g_0)^2 \int_0^1 |f_k|^2 dx
+  \int_0^1 |g_k |^2 dx + \frac{1}{2} \| D_t^{\alpha} u_k^n \|^2_{L_2(0,1)}.
\]
For the second summand in the right-hand side of \eqref{121} we have (see Eq. \eqref{14})
\[
J_2(t)\leq \frac{1}{4} \| D_t^{\alpha} u_k^n \|^2_{L_2(0,1)} + \frac{f_0^2}{2\varepsilon} \int_{0}^{1} |f_k|^2 \sum_{j=1}^{\infty} \lambda^{\frac{5}{2}+\varepsilon}_j |u_j^{n-1}|^2  dx.
\]
Therefore, inequality \eqref{121} can be rewritten as
\[
\frac{1}{4}\| D_t^{\alpha} u_k^n \|^2_{L_2(0,1)}+ \frac{1}{2} D_t^{\alpha} \| (u_k^n)_x \|^2_{L_2(0,1)} + \lambda_k \frac{1}{2} D_t^{\alpha} \| u_k^n \|^2_{L_2(0,1)}
\]
\[
\leq f_0 ^2(\psi_0 + g_0)^2 \|f_k\|^2_{L_2(0,1)}
+  \|g_k\|^2_{L_2(0,1)} +  \frac{f_0^2}{2\varepsilon} \int_{0}^{1} |f_k|^2 \sum_{j=1}^{\infty} \lambda^{\frac{5}{2}+\varepsilon}_j |u_j^{n-1}|^2  dx.
\]
Now we apply the operator $J_t^\alpha$ to both parts of inequality (\ref{15}) and use Lemma \ref{J}. Then
\[
\frac{1}{4}J_t^\alpha\| D_t^{\alpha} u_k^n \|^2_{L_2(0,1)}+ \frac{1}{2} \| (u_k^n)_x \|^2_{L_2(0,1)} + \lambda_k \frac{1}{2} \| u_k^n \|^2_{L_2(0,1)}
\]
\[
\leq \frac{1}{2} \| (\varphi_k)_x \|^2_{L_2(0,1)} + \lambda_k \frac{1}{2} \| \varphi_k \|^2_{L_2(0,1)}+f_0 ^2(\psi_0 + g_0)^2 J_t^\alpha\|f_k\|^2_{L_2(0,1)}
+  J_t^\alpha\|g_k\|^2_{L_2(0,1)}
\]
\[
+  \frac{f_0^2}{2\varepsilon} J_t^\alpha\int_{0}^{1} |f_k|^2 \sum_{j=1}^{\infty} \lambda^{\frac{5}{2}+\varepsilon}_j |u_j^{n-1}|^2  dx.
\]
First, we omit the second and third terms on the left-hand side of this inequality. Then, summing over \( k \) from 1 to \( \infty \) and applying the estimates (\ref{integral}), we obtain (\ref{Dalpha1}).
\end{proof}

\section{Convergence}
Let us show that the sequences \( u_k^n \), \( (u_k^n)_x \) and \( D_t^\alpha u_k^n \), \( n=1,2,..., \), are fundamental with respect to the corresponding norm for all \( k \geq 1 \).

Since the initial-boundary value problem \eqref{9}–\eqref{10} is linear with respect to \( u_k^n \), by repeating the reasoning above, we can establish estimates \eqref{18}, \eqref{u_x}, and \eqref{Dalpha} for \( u_k^{n+p} - u_k^n \), \( p = 1,2,3, \dots \). Indeed, for example,  write the equality (\ref{9}) for \(u_k^{n+1} \) and subtract (\ref{9}) from the resulting equality. Then we will have
\begin{equation}\label{Dalpha}
\int_0^1 D_t^\alpha (u_k^{n+1}-u_k^n)\, w \, dx + \int_0^1 (u_k^{n+1}-u_k^n)_x \,w_x \, dx + \lambda_k \int_0^1 (u_k^{n+1}-u_k^n) \,w \, dx 
\end{equation}
\[
=  - \int_0^1 \frac{f_k(t, x)}{f(t, x, l_0)} \sum_{j=1}^\infty \lambda_j (u_k^{n}-u_k^{n-1})(t, x) \sin j l_0 \,w \, dx,
\]
for all \( w \in \dot{W}_2^1(0, 1) \), with the initial condition:
\begin{equation}\label{101}
(u_k^{n+1}-u_k^n)(0, x) = 0.
\end{equation}
Now, repeating the same reasoning as in the proof of \eqref{17}, but using \eqref{101} instead of the initial condition \eqref{10}, we arrive at the inequality
\[
\sum_{k=1}^{\infty} \lambda_k^{\frac{5}{2} + \varepsilon} \| u_k^{n+1} - u_k^n \|^2_{L_2(0,1)} \leq\frac{1}{2}\sum_{k=1}^{\infty} \lambda_k^{\frac{5}{2} + \varepsilon} \| u_k^n - u_k^{n-1} \|^2_{L_2(0,1)}.
\]
Hence, keeping in mind that \( u_k^0 = 0 \) for all \( k \geq 1 \), from \eqref{u1}, we obtain
\[
\sum_{k=1}^{\infty} \lambda_k^{\frac{5}{2} + \varepsilon} \| u_k^2 - u_k^1 \|^2_{L_2(0,1)} \leq A_0. 
\]
For any \( n \geq 2 \) it follows from the last two estimates:
\[
\sum_{k=1}^{\infty} \lambda_k^{\frac{5}{2} + \varepsilon} \| u_k^{n+1} - u_k^n \|^2_{L_2(0,1)} \leq A_0 \left( \frac{1}{2} \right)^{n-2}.
\]
Hence, by writing $u_k^{n+p} - u_k^n=u_k^{n+p} - u_k^{n+p-1}+ u_k^{n+p-1}- u_k^{n+p-2}+\cdots$, we have
\begin{equation}\label{fundamental}
\sum_{k=1}^{\infty} \lambda_k^{\frac{5}{2} + \varepsilon} \| u_k^{n+p} - u_k^n \|^2_{L_2(0,1)} \leq A_0 \left( \frac{1}{2} \right)^{n-2} \left( 1 + \frac{1}{2} + \dots + \frac{1}{2^{p}} \right)\leq A_0 \left( \frac{1}{2} \right)^{n-1}.
\end{equation}
This implies that the sequence \(\{u_k^n\}\) is fundamental in \(L_2(0,1)\). Consequently, for each \(k \geq 1\), there exist functions \(u_k(t,x) \in L_2(0,1)\) such that, for every \(0 < t \leq T\), as \(n \to \infty\), the following conditions hold:
\begin{equation}\label{converU}
\begin{cases}
u_k^n(t, x) \to u_k (t, x)\quad \text{in } L_2(0,1), \quad k=1,2, 
\cdots,\\
\sum_{k=1}^{\infty} \lambda_k^{\frac{5}{2} + \varepsilon} \| u_k^n(t, \cdot) - u_k(t, \cdot) \|_{L_2(0,1)}^2 \to 0.
\end{cases}
\end{equation}

If we take $ u_k^{n+p}-u_k^n $ instead of $u_k^n $ in the estimate \eqref{151}, then it takes the form (in this case, it is necessary to repeat similar reasoning to that given in the proof of (\ref{fundamental}))
\[
\| u_k^{n+p}-u_k^n \|_{L_2(0, 1)}^2 + 2 J_t^\alpha\| (u_k^{n+p}-u_k^n )_x \|_{L_2(0, 1)}^2 + 2 \lambda_k J_t^\alpha\| u_k^{n+p}-u_k^n  \|_{L_2(0, 1)}^2 
\]
\[
\leq J_t^\alpha \| u_k^{n+p}-u_k^n  \|_{L_2(0, 1)}^2 +  J_t^\alpha c_{2,2}^k(t).
\]
After simple calculations, using the estimate (\ref{integral}), (\ref{u_x}) can be rewritten as
\[
\sum\limits_{k=1}^\infty\int\limits_0^t\int\limits_0^1\left[| (u_k^{n+p}-u_k^n)_x |^2 + \lambda_k | u_k^{n+p}-u_k^n |^2\right]dx d\tau \leq  T\,\max_{t\in [0,T]} \sum\limits_{k=1}^\infty \| u_k^{n+p}-u_k^n \|_{L_2(0, 1)}^2
\]
\[
+\frac{T\,  f_0^2}{2 \varepsilon} \max_{x\in [0,1]} \int\limits_0^\pi |f(t,x,y)|^2 dy\, \sum_{k=1}^{\infty} \lambda_k^{\frac{5}{2} + \varepsilon} \| u_k^{n+p-1} - u_k^{n-1} \|^2_{L_2(0,1)}, \quad t\in [0,T]. 
\]
The estimate (\ref{fundamental}) implies that the sequences on the left-hand side are fundamental. Hence, for every \(0 < t \leq T\), as \(n \to \infty\), we have
\begin{equation}\label{converUx}
\begin{cases}
\int\limits_0^t \| (u_k^n)_x(\tau, \cdot) - (u_k)_x(\tau, \cdot) \|_{L_2(0,1)}^2 d\tau \to 0, , \quad k=1,2, 
\cdots,\\
\int\limits_0^t \lambda_k\| u_k^n(\tau, \cdot) - u_k(\tau, \cdot) \|_{L_2(0,1)}^2 d\tau\to 0, , \quad k=1,2, 
\cdots,\\
\int\limits_0^t\sum_{k=1}^{\infty}  \| (u_k^n)_x(\tau, \cdot) - (u_k)_x(\tau, \cdot) \|_{L_2(0,1)}^2 d\tau \to 0, \\
\int\limits_0^t\sum_{k=1}^{\infty} \lambda_k \| u_k^n(\tau, \cdot) - u_k(\tau, \cdot) \|_{L_2(0,1)}^2 d\tau\to 0.
\end{cases}
\end{equation}

Let us take $w= D^\alpha_t(u_k^{n+1}- u_k^n)$ in equality (\ref{Dalpha}) and repeat the reasoning given in the proof of (\ref{Dalpha1}). Next, using similar arguments that were used in the proof of (\ref{fundamental}), we will verify that the sequence $ \sum\limits_{k=1}^\infty\int\limits_0^t\int\limits_0^1| D_t^\alpha u_k^n |^2dxd\tau$ is fundamental. Therefore, for every \(0 < t \leq T\), as \(n \to \infty\), one has
\begin{equation}\label{fundamental2}
\begin{cases}
\int\limits_0^t \| D_t^\alpha u_k^n (\tau, \cdot) - D_t^\alpha u_k(\tau, \cdot) \|_{L_2(0,1)}^2 d\tau \to 0, \quad k=1,2, 
\cdots, \\
    \int\limits_0^t\sum_{k=1}^{\infty}  \| D_t^\alpha u_k^n (\tau, \cdot) - D_t^\alpha u_k(\tau, \cdot) \|_{L_2(0,1)}^2 d\tau \to 0.
    \end{cases}
\end{equation}

 Let us prove convergence in \(L_2(0,1)\) of the right-hand side of the equation (\ref{9}). We have
\begin{equation}\label{right part}
\left|\sum\limits_{k=1}^{\infty}\lambda_k \big(u_k^{n} - u_k\big)\sin{kl_0}\right|^2\leq \sum_{k=1}^{\infty} \lambda_k^{-\frac{1}{2} - \varepsilon} \cdot \sum_{k=1}^{\infty} \lambda_k^{\frac{5}{2} + \varepsilon} |u^n_k (t, x)-u_k (t, x)|^2 
\end{equation}
\[
\leq \frac{1}{2 \varepsilon} \sum_{k=1}^{\infty} \lambda_k^{\frac{5}{2} + \varepsilon} |u^n_k (t, x)-u_k (t, x)|^2, \quad 0 \leq x \leq 1, \quad t \geq 0.
\]
Integrating over the interval $(0,1)$ and taking into account (\ref{converU}), we obtain the required convergence.

\section{Proof of Theorem \ref{theorem 1}}
If we pass to the limit $n\to \infty$ in (\ref{18}), taking into account the convergence (\ref{converU}) and Parseval's equality, we come to the conclusion that for each $t\in [0,T]$ the function
\[
u (t, x, y) = \sum\limits_{k=1}^{\infty} u_k (t, x) \sin k y,
\]
is well defined element of $L_2(\Omega)$ and (see Lemma \ref{u_k})
\begin{equation}\label{uestimate}
\max_{[0,T]}\| u(t, \cdot, \cdot) \|_2^2\leq 2 A_0.
\end{equation}
Similarly, passing to the limit $n\to \infty$ with respect to (\ref{u_x}), taking into account the convergence (\ref{converUx}) we will have
 \[
\int\limits_0^t\left[|| u_x(\tau. \cdot, \cdot) ||^2 +|| u_y(\tau. \cdot, \cdot) ||^2\right]d\tau \leq \frac{\Gamma(\alpha)T^{1-\alpha}}{2}||\varphi||^2_2+ 3 A_0 T +\frac{T}{2} \, A_1 
\]
\begin{equation}\label{u_xu_y}
+\frac{T\, A_0 f_0^2}{2 \varepsilon} \max_{t, x} \int\limits_0^\pi |f(t,x,y)|^2 dy , \quad t\in [0,T]. 
\end{equation}
The estimate (\ref{Dalpha1}) implies
\begin{equation}\label{DalphaEstimate}
\int\limits_0^t|| D_t^\alpha u(\tau, \cdot, \cdot)||^2d\tau \leq \frac{\Gamma(\alpha)T^{1-\alpha}}{2} \left(\| \varphi_x \|^2_2 +  \| \varphi_y \|^2_2\right)
\end{equation}
\[
+T\max_t\left[f_0 ^2(\psi_0 + g_0)^2 \|f(t, \cdot, \cdot)\|^2_2
+ \|g(t, \cdot, \cdot)\|^2_2\right]+\frac{A_0 f_0^2}{\varepsilon} \max_{t, x} \int\limits_0^\pi |f(t,x,y)|^2 dy, \quad t\in [0,T]. 
\]

 Let us now show that the function \( h (t, x) \), defined by the formula \eqref{7}, is an element of $L_2(0,1)$ for each $t\in [0,T]$.
 \begin{лем}\label{hL2} One has
 \[
\int_0^1|h(t,x)|^2dx \leq 4 f_0^2\left(\psi_0^2 +  \frac{A_0}{\varepsilon}\right)+ 2 g_0^2, \,\, t\in [0,T].
\]
\end{лем}
\begin{proof} By definition (\ref{7}) we get
  \begin{equation}\label{h1}
|h(t,x)|^2\leq 2 f_0^2 \left(\psi_{0} + \left|\sum\limits_{k=1}^{\infty} \lambda_k u_k (t, x) \sin k l_0\right| \right)^2 + 2 g_0^2.
\end{equation}
On the other hand, one has
\[
\left| \sum_{k=1}^{\infty} \lambda_k \, u_k (t, x) \sin k l_0 \right|^2
\leq \sum_{k=1}^{\infty} \lambda_k^{-\frac{1}{2} - \varepsilon} \cdot \sum_{k=1}^{\infty} \lambda_k^{\frac{5}{2} + \varepsilon} |u_k (t, x)|^2 
\]
\[
\leq \frac{1}{2 \varepsilon} \sum_{k=1}^{\infty} \lambda_k^{\frac{5}{2} + \varepsilon} |u_k (t, x)|^2, \quad 0 \leq x \leq 1, \quad t \geq 0.
\]
Let us integrate this inequality with respect to $x\in (0,1)$. Then, by virtue of estimate (\ref{18}), we will have
\begin{equation}\label{h2}
\int_0^1\left| \sum_{k=1}^{\infty} \lambda_k \, u_k (t, x) \sin k l_0 \right|^2dx \leq \frac{A_0}{\varepsilon}.
\end{equation}
Now we integrate (\ref{h1}) with respect to $x$, and apply the inequality $(a+b)^2 \leq 2 a^2+ 2b^2$. Then by (\ref{h2}), we get the required estimate.  
\end{proof}

We integrate \eqref{9} over $t\in [0,T]$, and then pass to the limit $n \rightarrow \infty$ and, taking into account \eqref{converU} - \eqref{right part}, we obtain the following equalities valid for arbitrary \( w \in \dot{W}_2^1(0, 1) \) and for all $k$:
 \begin{equation}\label{22}
\int_0^t\left[( D_t^{\alpha} u_k, w ) + ( (u_k)_x, w_x ) + \lambda_k ( u_k, w )\right] d\tau =\int_0^t\left[( g_k, w )+( f_kh, w )\right] d\tau.  
\end{equation}
Next, we take the function $v(x,y) \sin ky$, $v\in \dot{W}_2^1(\Omega)$ as $w$, integrate the equality (\ref{22}) over $y\in [0, \pi]$ and sum over $k$ from 1 to $\infty$. Then for all $t\in [0, T]$ one has
 \begin{equation}\label{solution_t}
\int_0^t \int_\Omega \left[D_t^\alpha u v +  u_x v_x + u_y v_y \right]\, dx \, dy d\tau = \int_0^t\int_\Omega \left[f h v + g v \right]\, dx \, dy d\tau.
\end{equation}
Note that here we have used the obvious equality $(u_{yy}, v)= (u_y, v_y)$.

Now we note the obvious fact that if $s(t)$ is integrable in any subset $(0,t)$ of $[0,T]$ and $\int_0^t s(\tau) d\tau =0$, then $s(t)=0$ is almost everywhere on $[0,T]$. Therefore, the equality (\ref{solution_t}) coincides with (\ref{equation}). In other words, the function $u(t,x,y)$ defined by formula (\ref{6}) is the weak solution of the inverse problem, i.e., it satisfies all conditions of Definition \ref{opr1} (see Remark \ref{zam.2}). The estimates \eqref{uestimate1} - (\ref{h}) are proved above (see estimates (\ref{uestimate}) - (\ref{DalphaEstimate}), Lemma \ref{hL2} and Remark \ref{zam.2}). The estimate (\ref{182}) is a consequence of inequality (\ref{181}).

To prove Theorem \ref{theorem 1}, it remains to prove the uniqueness of the solution of the inverse problem. 

Suppose that there are two solutions to the inverse problem: \( (u_1, h_1) \) and \( (u_2, h_2) \).  
Let us denote \( u = u_1 - u_2 \) and \( h = h_1 - h_2 \). Then the function \( u \) satisfies all the conditions of Definition \ref{opr1} with the functions
\( \varphi(x,y) \equiv 0 \), \( g(t,x,y) \equiv 0 \). Let us denote
\[
u_k(t, x) = \frac{2}{\pi} \int u(t, x, y) \sin k y \, dy.
\]
Consider the sum that is involved in defining the function $h(t,x)$ (see (\ref{7})):
\[
\sum_{k=1}^{\infty} \lambda_k u_k (t,x) \sin k l_0.
\]
By virtue of estimate (\ref{h2}), this function
exists almost everywhere in $(0,1)$ for all $t\in [0, T]$. This means that the function
\[
h(t,x) = - \frac{\sum\limits_{k=1}^{\infty} \lambda_k u_k (t,x) \sin k l_0}{f(t,x,l_0)}
\]
(see \eqref{7}) is correctly defined. Then the following estimate holds for the Fourier coefficients $u_k$ (see \eqref{12}):

\[
\frac{1}{2} D_t^\alpha \| u_k \|^2_{L_2(0,1)} + \| (u_k)_x \|^2_{L_2(0,1)} + \lambda_k \| u_k \|^2_{L_2(0,1)}\leq 
\]
\[
\leq \left| \int_0^1 u_k\frac{f_k (t,x)}{f(t,x,l_0)} \sum_{j=1}^{\infty} \lambda_j u_j \sin j l_0 \,dx \right|
\]
 Note that $M_k (t,x) \equiv 0$. Now, repeating the same reasoning as in the proof of \eqref{18}, we obtain:
\[
\sum_{k=1}^{\infty} \lambda_k^{\frac{5}{2} + \varepsilon} \| u_k \|^2_{L_2(0,1)} \leq 0
\]
(note that $A_0= 0$). Hence,
\[
u_k (t,x) = 0 \quad \text{a.e. in } (0,1) \text{ for all  } k \geq 1 \text{ and } t\in [0, T].
\]
 This means that $h(t,x) = 0$ \text{a.e. in} $(0,1)$ for all $t\in [0, T]$. By virtue of the completeness of the system $\{ \sin j y \}$, it follows that
\[
u(t,x,y) = 0 \quad \text{for almost all}\,\,  x\in (0,1), \, y \in (0,\pi) \text{ and } t\in [0, T].
\]
Thus, the theorem is completely proved.

\begin{зам}\label{zam.3} Note that the Theorem \ref{theorem 1} remains valid for parabolic equations as well. In this case, the usual equality $\frac{d}{dt} u^2= 2 u \frac{d}{dt} u$ should be used instead of Alikhanov's estimate.
\end{зам}

\section{Conclusions}
 
In this paper, we investigate a novel inverse problem for subdiffusion equations, specifically the problem of determining a source function that depends on both time and a subset of the spatial variables. To the best of our knowledge, such a problem for equations involving fractional derivatives has not been previously explored. We propose a method for solving this problem, drawing partially on the work of B.A. Bubnov. The existence and uniqueness of a weak solution to the inverse problem are established, along with coercive estimates.

The elliptic part of the equation is given by $u_{xx} + u_{yy}$, and the unknown function depends on a portion of the spatial variables, denoted as $h(t, x)$. The proposed method can be extended, without significant modifications, to the case where $x \in \mathbb{R}^n$ and $y \in \mathbb{R}^m$, and the equation incorporates two arbitrary elliptic operators acting separately with respect to the variables $x$ and $y$.

The choice to construct a weak solution is motivated by the relative ease of obtaining a priori estimates in this framework. Alternatively, if we consider the classical solution form of the problem (\ref{9}), (\ref{10}) (e.g.,  see \cite{AMux2}), it is possible to prove the existence of a classical solution to the inverse problem. However, this will be the focus of a future article.

\

\begin{center}
ACKNOWLEDGEMENTS    
\end{center}

 The authors acknowledge financial support from the Ministry of Innovative Development of the Republic of Uzbekistan, Grant No F-FA-2021-424.


\end{document}